\documentclass[12pt]{article}
\usepackage{amsmath,amssymb,amsthm}
\newtheorem{theorem}{Theorem}[section]

\newtheorem{cor}[theorem]{Corollary}

\theoremstyle{definition}

\theoremstyle{remark}

\numberwithin{equation}{section}

\newcommand{\myqed}{\hfill $\diamondsuit$   \medskip}

\newenvironment{proof1}{{\bf Proof :} }{\myqed}

\begin{document}

\title{{Stable Ranks, $K$-Groups 
and Witt Groups of some Banach and $C^{\ast}$-Algebras}
\thanks{Preprint, to appear in Contemporary Math.}}

\author{{C. BADEA}\\
$\; $\\
Math\'ematiques, UMR 8524 au CNRS\\
Universit\'e des Sciences et Technologies de Lille\\
F-59655 Villeneuve d'Ascq, France\\
E-mail : {\tt badea@gat.univ-lille1.fr}\\
URL : {\tt www-gat.univ-lille1.fr/\~{}badea}}
\date{ }
\maketitle
\thispagestyle{empty}
\begin{quote}
{\bf Abstract :} We show 
that certain dense and spectral invariant subalgebras of a $C^{\ast}$-algebra 
have the same bilateral Bass stable rank. This is a partial answer for 
(a version of) an open problem raised by R.G. Swan. 
Then, for certain Banach algebras, 
we indicate when the homotopy 
groups $\pi_{i}(GL_{n}(A))$ stabilize for large $n$. This is an 
improvement of a result due to G.~Corach and A.~Larotonda. 
Using some results due 
to M.~Karoubi, we show the isomorphism of the Witt group 
of a symmetric 
Banach algebra with the $K_{0}$-group of its 
enveloping $C^{\ast}$-algebra. The question if this is true for all 
involutive Banach algebras was raised by 
A.~Connes.
\end{quote}
\begin{quote}
{\bf R\'esum\'e :} On d\'emontre que certaines sous-alg\`ebres 
denses et pleines d'une $C^{\ast}$-alg\`ebre ont le m\^eme rang stable (de Bass)
bilat\`ere. Ceci est une r\'eponse partielle \`a un probl\`eme 
ouvert de R.G. Swan. Pour certaines alg\`ebres de Banach $A$,
les groupes d'homotopie $\pi_{i}(GL_{n}(A))$ sont 
stables pour $n$ assez grand.
Ceci est une am\'elioration d'un r\'esultat de 
G. Corach et A. Larotonda. En utilisant des r\'esultats 
dus \`a M. Karoubi, on d\'emontre  
l'isomorphisme du groupe de 
Witt d'une alg\`ebre de Banach sym\'etrique avec le groupe $K_{0}$ de 
son $C^{\ast}$-alg\`ebre enveloppante. La question de savoir si cet 
isomorphisme a lieu pour toutes les alg\`ebres de Banach 
involutives \`a \'et\'e soulev\'ee par A. Connes.
\end{quote}
{\bf 1991 Subject Classification} : Primary 46H05, 46H10

\newpage

\tableofcontents
\section{Introduction}
 There are several notions of stable ranks for topological algebras. 
 We will discuss here some of them for some Banach and 
 $C^{\ast}$-algebras and their incidence in the study of certain 
 groups related to $K$-Theory. Their definitions and main properties
 are surveyed in the next 
 section. 
 
The problems we are dealing 
with are 
\begin{itemize}
	\item  Swan's problem for subalgebras of $C^{\ast}$-algebras and for the 
	bilateral Bass stable rank ;

	\item  stabilization 
of the homotopy groups of the general 
linear group of a Banach algebra ;

	\item  the isomorphism between the 
Witt group and the $K_{0}$-group of the enveloping $C^{\ast}$-algebra 
of a given symmetric Banach algebra.
\end{itemize}

We present now the motivation for these problems and the results we 
will prove.

\subsection{Swan's problem for the bilateral Bass stable rank}
The density theorem in $K$-Theory implies that a Banach  dense and 
spectral invariant subalgebra $A$ 
of a Banach algebra $B$ has the same $K$-theory as $B$, that is the 
inclusion morphism $j : A \to B$ induces isomorphisms $j_{\ast} : 
K_{i}(A) \to K_{i}(B)$, $i = 0,1$. Recall that $A$ is called \emph{spectral 
invariant} in $B$ if $a \in A$ with $a$ invertible in $B$ imply that 
$a$ is invertible in $A$. Spectral invariant Banach subalgebras are 
closed under the holomorphic functional calculus, that is, for $a \in 
A$ and $g$ holomorphic in a neighborhood of the spectrum of $a$ in 
$B$, the element $g(a)$ of $B$ lies in $A$. 

It will be recalled in the next section that stable ranks stabilize 
the $K$-groups. The question whether dense spectral invariant 
subalgebras have the same stable rank arises. It was Richard G. Swan 
\cite[p.206]{Sw} who raised this question for the Bass stable rank 
and for the projective stable rank. In this note we consider Swan's 
problem for a variant of the Bass stable rank, called bilateral Bass 
stable rank \cite{Bad}. 

Swan's problem for several stable ranks was considered in \cite{Bad}. 
A consequence of the main result there gives a positive answer of 
Swan's problem for the Bass stable rank in the case when $B$ is a 
$C^{\ast}$-algebra and $A$ is a dense and 
spectral invariant $\ast$-subalgebra, which is a Fr\'echet 
$Q$-algebra (cf. \cite{Bad} for definitions) in its own topology and 
closed under $C^{\infty}$-functional calculus of selfadjoint elements.

We give in this note a sufficient condition for $\ast$-subalgebras 
of $C^{\ast}$-algebras to have the same bilateral Bass stable rank. 
In particular, we prove that dense spectral invariant $(D_{p})$ and 
$\ast$-subalgebras of $C^{\ast}$-algebras have the same bilateral 
Bass stable rank. The notion of $(D_{p})$-subalgebras of Banach 
algebras was recently considered by Kissin and Shulman \cite{KiSh2}.
 
\subsection{Computing homotopy groups}
Stability theorems for the $K$-groups
$K_{0}(A)$ and $K_{1}(A)$ 
in terms of the Bass stable rank are stated in the next section. Moreover, 
if the Bass stable rank of $A$ is finite, the homotopy groups 
$\pi_{i}(GL_{n}(A))$ stabilize for large $n$.  Here $GL_n(A)$ is the 
group of invertible
elements of $M_n(A)$, the set of all 
$n\times n$ matrices with entries in $A$. 

To be more specific, it was 
proved by G.~Corach and A.R.~Larotonda \cite{CoLa2} 
that the map 
$$\pi_i(GL_{n-1}(A)) \rightarrow \pi_i(GL_{n}(A))$$ 
between
the homotopy groups is surjective for $n \geq Bsr(A) + i+1$ 
and injective for $n \geq
Bsr(A) + i + 2$. For other results of this type we refer to \cite{Rie87},
\cite{Tho91}, \cite{Srd84}, \cite{Srd94}, \cite{Zha91}, \cite{Zha93}. 

We prove here a
stability result in terms of the connected stable rank and 
the general stable rank
which is a slight improvement of the result of Corach and Larotonda. 

The following consequence for commutative Banach algebra is obtained. 
Suppose that $A$ is unital and commutative. Then the map 
$$\pi_i(GL_{n-1}(A)) \rightarrow \pi_i(GL_{n}(A))$$
is surjective for $n \geq Bsr(A) + [i/2] + 2$ 
and injective for $n \geq
Bsr(A) + [(i+1)/2] + 2$.

\subsection{Computing the Witt groups}
We will consider in this section symmetric 
 Banach algebras, commutative or not. Recall that an involutive Banach 
 algebra is said to be \emph{symmetric} if every element of the form 
 $x^{\ast}x$ has spectrum included in positive real closed half-line. 
 The symmetry is equivalent in the commutative case to 
 $m(x^{\ast}) = \overline{m(x)}$ for all characters $m$. 

The following problem has been raised by A.~Connes in his book 
\cite{Co}.

It is known (see \cite{Co}) that for $C^{\ast}$-algebras $A$ one has 
an isomorphism between the Witt group $W_{0}(A)$ and the 
$K_{0}$-group $K_{0}(A)$ 
of $K$-theory. We refer to the next section for definitions. 
In \cite{Co} it is asked if, for an 
involutive Banach algebra $A$, the Witt group $W_{0}(A)$ is 
isomorphic to $K_{0}(C^{\ast}(A))$, the $K_{0}$-group of the 
enveloping $C^{\ast}$-algebra of $A$. J.-B.~Bost (cf. \cite{Co}) 
proved this for \emph{commutative} involutive $A$. We show in the 
present note that an affirmative answer (even for higher Witt groups) 
follows for symmetric Banach involutive algebras 
from the work of M.~Karoubi \cite{Kar2}. 

\subsection{Acknowledgments}
The present note is a modified version of my talk at the Third Conference on 
Function Spaces in Edwardsville, Illinois (May 1998). I wish 
to thank K. Jarosz 
for his invitation and the referee for a careful reading.

\section{Background}
A Banach algebra will mean a complex Banach algebra. The unity is 
always denoted by $e$. 

\subsection{Stable ranks}
For a Banach algebra $A$ we denote by $Lg_{n}(A)$ the set of all 
$n$-tuples $(a_{1}, \ldots, a_{n}) \in A^n$ with the property 
$Aa_{1} + \cdots Aa_{n} = A$, i.e. generating $A$ as a left ideal. By 
definition, the (left) \emph{Bass stable rank} $Bsr(A)$ of $A$ is the 
smallest positive integer $n$ such that the following condition 
\begin{quote}
$(Bsr)_n$ \quad for every $(a_1, \ldots ,a_{n+1})\in Lg_{n+1}(A)$,
there exists
$(c_1,\ldots ,c_n) \in A^n$ such that $(a_1+c_1a_{n+1},\ldots
,a_n+c_na_{n+1})\in Lg_n(A)$, 
\end{quote}
holds, or infinity if no such number $n$ exists. 
Note that in some of the earlier papers
the notion of the Bass stable rank was defined using different
indexing conventions. The condition
$(Bsr)_n$ was devised by H.~Bass in order to 
determine values of $n$ for which every
matrix in  $GL_n(A)$ can be   row reduced by addition operations with 
coefficients
from $A$ to a matrix with the same last row and column as the identity matrix and to
obtain stability results in $K$-theory.

It can
be proved that $(Bsr)_n$ implies $(Bsr)_{n+1}$ (\cite{Vas},\cite{Kru}) 
and that
the Bass stable rank of $A$ equals the stable rank of the opposed algebra 
$A^\circ$
(\cite{Vas},\cite{War}). That is, the left Bass stable rank defined
above is equal to the right Bass stable rank which can be defined in a 
similar way
using the set $Rg_n(A)$ of all $n$-tuples $(a_1,\ldots ,a_n)$
generating $A$ as a right ideal. If $A$ is non-unital, we define the 
Bass stable rank of $A$ as the Bass stable rank of the algebra 
$A_{+}$ obtained by adjoining a unit element to $A$. 

A variant of the Bass stable rank was introduced in \cite{Bad}. We 
call \emph{bilateral Bass stable rank} $bBsr(A)$ of the unital Banach 
algebra $A$ the 
smallest positive integer $n$ such that the following condition holds 
for $k \geq n$
\begin{quote}
$(bBsr)_k$ \quad for every $(a_1, \ldots ,a_{k+1})\in Lg_{k+1}(A)$,
there exist
$(c_1,\ldots ,c_k) \in A^k$, $(d_1,\ldots ,d_k) \in A^k$ 
such that $(a_1+c_1a_{k+1}d_{1},\ldots
,a_k+c_ka_{k+1}d_{k})\in Lg_k(A)$, 
\end{quote}
or infinity if no such number exists. We have $bBsr(A) \leq Bsr(A)$, 
with equality for commutative $A$. There are \cite{Bad} $C^\ast$-algebras $A$ 
with $bBsr(A) \neq Bsr(A)$. 

M.A. Rieffel \cite{Rie} introduced the notion of topological stable rank 
as follows~:
the (left) \emph{topological stable rank} $tsr(A)$ of $A$ is the
smallest positive integer $n$ such that $Lg_n(A)$ is dense in $A^n$,
or infinity if no such number exists. If
$Lg_n(A)$ is dense in $A^n$, then $Lg_m(A)$ is dense in $A^m$ for every $m\geq n$.
A symmetric notion, the right topological stable rank $rtsr$ can be defined 
by considering the set
$Rg_n(A)$ instead of $Lg_n(A)$. The left and the right topological stable ranks
coincide for Banach algebras with a continuous involution. 
It is an open question
\cite[Question 1.5]{Rie} if $tsr(A)$ equals $rtsr(A)$ for 
all Banach algebras $A$. If $A$ is non-unital, we define the 
topological stable rank of $A$ as the topological stable rank of 
$A_{+}$.
We refer to \cite{Rie} for several properties of $tsr$. 

\subsection{Other stable ranks}
We mention briefly other notions of stable ranks. The {\it{connected stable
rank}\/} $csr(A)$ of the Banach algebra $A$ is \cite{Rie} the least integer 
$n$ such
that $GL_k(A)_0$ acts transitively (by left multiplication) on $Lg_k(A)$ for 
every  $k\geq n$, or, equivalently \cite{Rie}, the least integer $n$ such that
$Lg_k(A)$ is connected for every $k\geq n$. We put $csr(A) =
\infty$ if no such $n$ exists. This notion is left-right symmetric 
\cite{CoLa3}.
Recall that $GL_k(A)_0$ is the connected 
component of $GL_k(A)$ containing the identity.

The {\it{left (right) general stable rank}\/} of $A$ is defined 
\cite{Rie} as
the smallest integer $n$ such that $GL_k(A)$ acts on the left 
(right) transitively
on $Lg_k(A)$ for all $k\geq n$. If no such integer exists 
we set $gsr(A) = \infty$. For $C^\ast$-algebras it is related
to the
cancellation property for finitely generated projective $A$-modules. 
For
instance, the right general stable rank is the smallest positive integer 
$n$ such that
$W\oplus A \cong A^k$ for some $k\geq n$ implies $W\cong A^{k-1}$, 
whenever $W$ is a
finitely generated projective left $A$-module. 
This notion is left-right symmetric
\cite{CoLa3} and we will denote by $gsr(A)$ the common value. 

\subsection{Properties of stable ranks}
It was proved in \cite{Bad} that condition $(bBsr)_{n}$ holds for $A$ 
if and only if every onto unital algebra morphism $f : A \to B$, $B$ a 
Banach algebra, induces an onto mapping $f_{n} : Lg_{n}(A) \to 
Lg_{n}(B)$. A similar characterization holds for the (left) Bass stable rank 
by replacing onto algebra morphisms from $A$ with onto module 
morphisms of left $A$-modules from ${}_{A}A$, which is $A$ 
viewed as a left $A$-module. 
Also \cite{Bad}, $tsr(A) \leq n$ if and only if 
for every $\varepsilon > 0$ and every
$(a_1,\ldots ,a_{n+1})\in Lg_{n+1}(A)$, there exists $(c_1,\ldots
,c_n)\in A^n$ such that 
$(a_1+c_1a_{n+1},\ldots ,a_n+c_na_{n+1})\in Lg_n(A)$ and
$\| c_ia_{n+1}\| \leq \varepsilon$ for
$i=1,\ldots ,n$. This yields 
$bBsr(A) \leq Bsr(A) \leq tsr(A)$ for all Banach algebras $A$. 
For $C^\ast$-algebras 
we have
$Bsr(A) = tsr(A)$ as was shown by Herman and Vaserstein \cite{HeVa}.
For the unital commutative $C^\ast$-algebra $C(X)$ we have
$$bBsr(C(X)) = Bsr(C(X)) = tsr(C(X)) = [(\dim X)/2] + 1,$$
where $\dim X$ is the \v{C}ech-Lebesgue 
covering dimension of $X$ \cite{Pea}.

We have 
$$gsr (A) \leq csr (A) \leq 1 + Bsr (A) \leq 1 + tsr (A) \; .$$

\subsection{Stable ranks and $K$-theory}
The Bass and topological stable ranks are useful for
stability results in the $K$-theory of topological algebras. We state here some of
them in terms of $Bsr$~; since $Bsr(A) \leq tsr(A)$, we can replace in these
statements the Bass stable rank with the topological one. 
It was proved by G. Corach
and A.R. Larotonda \cite{CoLa2} that, for a Banach algebra $A$ with unit $e$,
the map
$$
GL_n(A)/GL_n(A)_0\ni a\mapsto 
\left(
\begin{array}{cc} 
a & 0\\ 
0 & e
\end{array}
\right)
\in GL_{n+1}(A)/GL_{n+1}(A)_0$$
is bijective for $n\geq 1 + Bsr(A)$. Thus the topological
$K_1$-group $K_1(A)$ of $A$, defined as the direct limit of
$GL_n(A)/GL_n(A)_0$ under these inclusions, stabilizes if $Bsr(A)$
is finite. A similar result can be stated for the (topological)
$K_0$ group of a Banach algebra. Indeed, $K_0(A)$ can be written
\cite{Kar1} as the direct limit of $\bar P_{2n}(A)$ under the inclusions
$$p\mapsto
\left(
\begin{array}{ccc}
p & 0 & 0\\
0 & e & 0\\
0 & 0 & 0
\end{array}
\right) \; ,
$$
where $\bar P_{2n}(A)$ is the quotient of $P_{2n}(A)$, the subset of indempotents
of $M_{2n}(A)$, given by the equivalence relation
$$p_1 \approx p_2 \Leftrightarrow \exists \alpha\in GL_{2n}(A),
\quad \alpha p_1\alpha^{-1} = p_2\;.$$
Then the map $\bar P_{2n}(A)\rightarrow\bar P_{2n+2}(A)$ is bijective for 
$n\geq Bsr(A)$ \cite{Cor86}.

\subsection{Witt groups}
Let $A$ be an involutive unital Banach algebra. In 
Hermitian $K$-theory we denote by $L_{0}(A)$ the Grothendieck 
group of  finitely generated, projective $A$-modules endowed with 
non-degenerate quadratic forms. In order to define higher $L$-groups 
$L_{n}(A)$, we consider the orthogonal group $O_{k,k}(A)$ of isometries 
of $A^{2k}$ endowed with the standard hyperbolic form. It can be 
described as the group of $2k\times 2k$ matrices $X$ of the form
$$X = 
\left( \begin{array}{cc}
	M & N  \\
	P & Q
\end{array} \right) \; ,$$ 
where $M, N, P, Q \in M_{n}(A)$ and $X'X = XX' = I_{2n}$. Here $X'$ is
the following matrix 
$$\left( \begin{array}{cc}
	{}^tQ^{\ast} & {}^tN^{\ast}  \\
	{}^tP^{\ast} & {}^tM^{\ast}
\end{array} \right),$$
while $I_{2n}$ is the identity matrix. 

Let $O(A)$ be the inductive limit of $O_{k,k}(A)$ under the 
inclusions 
$$a \to \left( \begin{array}{cc}
	a & 0  \\
	0 & e
\end{array} \right) \; .$$
The topological Hermitian $K$-groups $L_{n}(A)$ are defined as 
homotopy groups of $O(A)$, namely
$$L_{n}(A) = \pi_{n-1}(O(A))\; .$$
We refer to \cite{Kar2} and the references therein for further 
information. In particular, the hyperbolic functor induces a 
homomorphism $K_{n}(A) \to L_{n}(A)$. 

The $n$-th Witt group $W_{n}(A)$ is defined as 
$$W_{n}(A) = \textrm{Coker} (K_{n}(A) \to L_{n}(A))\; .$$
For $n = 0$ we obtain the classical Witt group $W_{0}$. 

\section{Computing $bBsr$ for dense subalgebras of $C^{\ast}$-algebras}
The following is the announced partial result for Swan's problem for 
the bilateral Bass stable rank.

\begin{theorem}\label{th:3.1}
	Let $(B,\|\cdot\|_{0})$ be a unital $C^{\ast}$-algebra and let 
	$(A,\|\cdot\|_{0})$ be a dense and spectral invariant $\ast$-subalgebra 
	containing the unity of $B$. 
	Let $p$ be a fixed positive integer. 
	Suppose that there is a family  $\{\|\cdot\|_{i}\}_{1\leq i\leq p}$ 
	of norms on $A$ and positive constants $C_{i}$, $1\leq i\leq p$, 
	such that
	$$\| x^2\|_{i} \leq C_{i}\|x\|_{i}\|x\|_{i-1} \; , \; 1\leq i\leq p,$$
	for any selfadjoint $x = x^{\ast}\in A$. 
	Suppose also that $(A,\|\cdot\|_{p})$ is a Banach algebra. Then 
	$bBsr(A) = bBsr(B)$.
\end{theorem}

\begin{proof1}
We have \cite[Theorem 4.7(i)]{Bad} $bBsr(A) \leq bBsr(B)$. For the 
reverse inequality, suppose $bBsr(A) = n$ and 
let $\mathbf{b} = (b_{1}, \ldots , b_{n+1}) \in Lg_{n+1}(B)$. We will prove 
that $\mathbf{b}$ is \emph{bilateral reducible} in $B$, that is there 
exist $(c_{1}, \ldots , c_{n})$, $(d_{1}, \ldots , d_{n})$ in $B^n$ 
such that 
$$(b_1+c_1b_{n+1}d_{1},\ldots
,b_n+c_nb_{n+1}d_{n})\in Lg_n(B) .$$
 
Let $J$ be the closed two-sided ideal in $B$ generated by $b_{n+1}$. 
It was proved in \cite[Theorem 4.7(ii)]{Bad} that $\mathbf{b}$ is 
bilateral reducible in $B$ if $J \cap A$ is dense 
in $J$. We prove now, using several techniques borrowed from 
\cite{KiSh1,KiSh2}, that $\overline{J \cap A} = J$ holds 
for all two-sided closed ideals 
of $B$.

Let $x = x^{\ast}$ be a selfadjoint element of A. Let $t$ be a real number. 
Let $k$ be the positive integer such that 
$2^{k-1} < |t| \leq 2^k$. Set $y = itx/2^k$. Let 
$$M = \max_{1\leq j\leq p} \{ C_{j}, \| \exp (y)\|_{j}\}.$$
Then $\| \exp (y)\|_{0} = 1$ and 
$$\| \exp (y)\|_{j} = \|\exp (itx/2^k)\|_{j} \leq \exp 
(|t|\|x\|_{j}/2^k) \leq M$$
for $j = 1, \ldots , p$. Using \cite[Lemma 1]{KiSh2} we get 
\begin{eqnarray*}
	\|\exp (itx)\|_{p} & = & \|\exp (2^ky)\|_{p} = \|\exp (y)^{2^k}\|_{p}  \\
	 & \leq & \|\exp (y)\|_{0}^{2^k-S(k,p)}\prod_{j=0}^{p-1}\|\exp 
	 (y)\|_{p-j}^{a(k,j)}C_{p-j}^{a(k,j+1)}  \\
	 & = & \prod_{j=0}^{p-1}\|\exp 
	 (y)\|_{p-j}^{a(k,j)}C_{p-j}^{a(k,j+1)} \leq M^b , 
\end{eqnarray*}
where 
$$a(k,j) = {k \choose j} \quad \mbox{and} \quad b = \sum_{j=0}^{p-1} 
\left[ a(k,j) + a(k,j+1)\right] \; .$$
For $k > 2p$ (and so for all $t$ such that $|t|$ is sufficiently 
large) we have (cf.~\cite[p. 416]{KiSh2})~:
$$\|\exp (itx)\|_{p} \leq M|2t|^{\psi(t)},$$
where 
$$\psi (t) = \frac{(\log_{2}|2t|)^{p-1}\log_{2}M}{(p-1)!}.$$
This implies that
$$\int_{-\infty}^{+\infty} \frac{\log \|\exp (itx)\|_{p}}{1+t^2}\, 
dt < +\infty.$$
It follows \cite{KiSh2} from Shilov's \cite[\S 15.6]{Na} 
condition of regularity 
that $A$ is locally normal in $B$. This means \cite{KiSh1} 
that there is a commutative 
Banach $\ast$-subalgebra $B(x)$ in $B$ such that $e$ and $x$ belong 
to $B(x)$ and such that $A(x) = A \cap B(x)$ is a dense normal 
subalgebra of $B(x)$. It was proved in \cite[Theorem 13]{KiSh1} that 
if $J$ is a closed two-sided ideal in $B$, then $J \cap A$ is dense 
in $J$. This completes the proof.
\end{proof1}

The condition 
$$\| x^2\|_{i} \leq C_{i}\|x\|_{i}\|x\|_{i-1} \; , \; 1\leq i\leq p,$$
for selfadjoint $x$ in the above theorem is a particular case of the 
condition $(D_{p})$ studied in \cite{KiSh2}. They introduced 
$(D_{p})$-subalgebras $A$ of Banach algebras $(B,\|\cdot\|_{0})$ as 
dense subalgebras of $B$ for which there exist norms 
$\{\|\cdot\|_{i}\}_{1\leq i\leq p}$ and positive constants 
$D_{i}$, $1\leq i\leq p$, such that $(B,\|\cdot\|_{p})$ is a Banach 
algebra and
$$(D_{p}) \quad \| xy\|_{i} \leq D_{i}\left( \|x\|_{i}\|y\|_{i-1} + 
\|x\|_{i-1}\|y\|_{i}\right), \; x,y \in B \; 1\leq i\leq p .$$
Thus the condition in Theorem \ref{th:3.1} is obtained from 
condition $(D_{p})$ for $y=x$ with $C_{j} = 2D_{j}$. We obtain the 
following consequence.

\begin{cor}
	Dense spectral invariant $(D_{p})$ and $\ast$-subalgebras of 
	$C^{\ast}$-algebras have same bilateral Bass stable rank.
\end{cor}

We refer to \cite{KiSh2} for several examples of $D_{p}$-subalgebras 
of $C^{\ast}$-algebras. For instance, the differential subalgebras of 
order $p$ studied by Blackadar and Cuntz \cite{BlCu} are 
$(D_{p})$-subalgebras and, for $p=1$, these classes coincide.

The above Corollary can be viewed as a noncommutative generalization 
of a result due to Vaserstein \cite{Vas}. He proved that the Banach 
algebras $C^k(X)$ of $k$-differentiable functions on a compact 
manifold $X$ have the same Bass stable rank as $C(X)$. It is 
well-known that $C^{\ast}$-algebras are noncommutative analogies of 
the algebras of continuous functions. $(D_{p})$-subalgebras of 
$C^{\ast}$-algebras can be 
viewed \cite{KiSh2} as noncommutative analogies of 
algebras of smooth functions. 
Note also that for commutative Banach algebras the bilateral Bass 
stable rank coincide with the usual Bass stable rank.

\section{Computing homotopy groups}
Denote by $\mathbf{S}^i$ the $i$-dimensional sphere in 
$\mathbb{R}^{i+1}$. 

The proof of the following stabilization theorem is based upon 
\cite{CoLa2} and \cite{Rie87}.

\begin{theorem}\label{th:homotopy} 
Let $A$ be a unital Banach algebra and let $n \geq csr(A)$. The
canonical morphism from $\pi_i(GL_{n-1}(A))$ into $\pi_i(GL_n(A))$ is
\begin{itemize} 
	\item[a)]	surjective for $n \geq csr(C(\mathbf{S}^i,A))$
	\item[b)]	injective for $n \geq gsr(C(\mathbf{S}^{i+1},A))$. 
\end{itemize}
In particular (cf. \cite{CoLa2}), the canonical morphism is 
\begin{itemize} 
	\item[a)]	surjective for $n \geq Bsr(A) + i+1$
	\item[b)]	injective for $n \geq
Bsr(A) + i + 2$. 
\end{itemize}
Moreover, if $A$ is commutative, then the canonical 
morphism from $\pi_i(GL_{n-1}(A))$ into $\pi_i(GL_n(A))$ is
\begin{itemize} 
	\item[a)]	surjective for $n \geq Bsr(A) + [i/2] +2$
	\item[b)]	injective for $n \geq Bsr(A) + [(i+1)/2] +2$.
\end{itemize}
 
 Here $[\cdot ]$ is the integer part.
\end{theorem}

\begin{proof1} 
By \cite{CoLa3}, the map $T : GL_n(A) 
\rightarrow Lg_n(A)$, $T(X) = Xe_n$, ($e_{n} = (0, \ldots , 0,1)$) 
is a Serre fibration and, 
with the same proof,
the map $T' : GL_n(A) \rightarrow Lc_n(A)$, $T'(X) = Xe_n$, is also a 
Serre fibration. Here $Lc_{n}(A)$ denotes the space of last columns 
of invertible $n\times n$ matrices with entries in $A$~: 
$$Lc_{n}(A) = \{ Me_{n} : M \in GL_{n}(A)\} \; .$$
The stability subgroup of $e_n$ induced by the action of $GL_n(A)$ 
under $Lg_n(A)$ 
and
$Lc_n(A)$ consists of matrices whose last column is $e_n$, 
that is matrices of the
form
$$
\left(
\begin{array}{cc}
x & 0 \\
c & e
\end{array}
\right) \; ,
$$
where $x \in GL_{n-1}(A)$ and $c$ is an arbitrary row in $A$ of length $n-1$. 
Viewing $GL_{n-1}(A)$ as a subset of $GL_n(A)$ via the embedding 
$$a \to  \left(
\begin{array}{cc}
a & 0 \\
0 & e
\end{array}
\right) \; ,
$$
we obtain a deformation retract of the stability subgroup of $e_n$ onto
$GL_{n-1}(A)$ by carrying the off-diagonal entry $c$ linearly to zero. Since $T'$ is a
Serre fibration, one has the homotopy exact sequence \cite[Ch. 7]{Spa66}
$$
\rightarrow \pi_{i+1}(Lc_n(A)) \rightarrow \pi_i(GL_{n-1}(A)) \rightarrow 
\pi_i(GL_n(A)) \rightarrow \pi_i(Lc_n(A)) \rightarrow \; .
$$
This long exact sequence ends with \cite{Spa66}
$$
\pi_0(GL_{n-1}(A)) \rightarrow \pi_0(GL_n(A)) \rightarrow \pi_0(Lc_n(A))
$$
viewed as pointed sets. The base points in the groups are taken to be their 
identity elements, while the base points in $Lg_n(A)$ and $Lc_n(A)$ are 
taken to be
$e_n$.    

For $n \geq csr(A)$ one has \cite{Bad} $n \geq gsr(A)$ and $n \geq Lccsr(A)$. 
Here $Lccsr(A)$, the \emph{last columns connected stable rank} is, 
by definition \cite{Bad}, the least integer $k$ such that 
for all $ n \geq k$ the 
set $Lc_{n}(A)$ is connected. Since \cite{Bad} 
$$Lc_{n}(A)_{0} = Lg_{n}(A)_{0} = GL_{n}(A)_{0}e_{n}\; ,$$
we have
$Lc_n(A) = Lg_n(A)$ and $\pi_0(Lc_n(A))$ is trivial. From the proof of 
\cite[Theorem 6.2]{CoLa2} it follows that there is a bijection 
$$\pi_i(Lg_n(A)) \to \pi_0(Lg_n(C(\mathbf{S}^i,A)))\; .$$ 
Since $n \geq
csr(C(\mathbf{S}^i,A))$, we get that $\pi_i(Lc_n(A)) = \pi_i(Lg_n(A))$ is also
trivial. From the long exact sequence above we obtain the surjectivity of
$$\pi_i(GL_{n-1}(A)) \rightarrow \pi_i(GL_n(A))\; .$$ 

It follows again from the homotopy exact sequence that the 
map 
$$\pi_i(GL_{n-1}(A)) \rightarrow \pi_i(GL_n(A))$$ 
is injective if and 
only if the map
$$\pi_{i+1}(GL_n(A)) \rightarrow \pi_{i+1}(Lc_n(A))$$ 
is surjective. Suppose that 
$n \geq gsr(C(\mathbf{S}^{i+1},A))$ and thus the map 
$$GL_n(C(\mathbf{S}^{i+1},A)) \to Lg_n(C(\mathbf{S}^{i+1},A))$$ 
is surjective.
Let $f : \mathbf{S}^{i+1} \rightarrow Lg_n(A)$ be a map preserving the base point
$f(1) = e_n$, and so representing an element of $\pi_{i+1}(Lg_n(A))$. Then $f$
can be identified with an element of $Lg_n(C(\mathbf{S}^{i+1},A))$ such that
$f(1) = e_n$. By the surjectivity, there is $g \in GL_n(C(\mathbf{S}^{i+1},A))$
such that $g(t)e_n = f(t)$, $t \in \mathbf{S}^{i+1}$. In particular, we have
$g(1)e_n = e_n$ and thus $g(1)^{-1}e_n = e_n$. Then the map 
$$h :
\mathbf{S}^{i+1}\ni t \rightarrow h(t) =
g(t)g(1)^{-1} \in GL_n(C(\mathbf{S}^{i+1},A))$$ 
is such that $h(t)e_n = f(t)$ and $h(1) = e$. Since $Lg_n(A) =
Lc_n(A)$, we obtain that the map 
$$\pi_{i+1}(GL_n(A)) \rightarrow \pi_{i+1}(Lc_n(A))$$ 
is surjective and thus 
$$\pi_i(GL_{n-1}(A)) \rightarrow \pi_i(GL_n(A))$$ 
is injective.

To obtain the first consequence (which is the result of Corach and 
Larotonda), note that 
$gsr(A) \leq csr(A) \leq 1 + Bsr(A)$ and \cite{CoLa2}
$Bsr(C(\mathbf{S}^i,A)) \leq i + Bsr(A)$. 

For the second consequence, we use a result due to F. D. Suarez 
\cite{SuaAlg} : if $A$ is commutative, then 
$$Bsr(C(\mathbf{S}^i,A)) \leq [i/2] + 1 + Bsr(A) \; .$$ 
This completes the proof.
\end{proof1}

It was conjectured that inequality $Bsr(C(\mathbf{S}^2,A)) \leq 
1 + Bsr(A)$ might be true \cite{Rie}, \cite{SuaAlg}, at least for 
$C^{\ast}$-algebras or for commutative Banach algebras. Suppose it holds 
for every complex Banach algebra. Then an inductive argument implies 
that $Bsr(C(\mathbf{S}^i,A)) \leq [(i+1)/2] + Bsr(A)$ for every $i$. 
This would lead 
to an improvement of Theorem \ref{th:homotopy}.

\section{Computing higher Witt groups}
The aim of this section is to show how some results due to M. 
Karoubi \cite{Kar2} yield an answer to the above-mentioned question of 
A. Connes in the case of symmetric Banach algebras. In fact, the 
following statement for higher Witt groups is true.

\begin{theorem}
	Let $A$ be a symmetric Banach $\ast$-algebra. For all $n \geq 0$, 
	$L_{n}(A)$ is isomorphic in a natural way to $K_{n}(A)\oplus 
	K_{n}(A)$. Therefore
	$$W_{n}(A) \simeq K_{n}(A) \simeq K_{n}(C^{\ast}(A))\; .$$
\end{theorem}

\begin{proof1}
The first statement was proved by M. Karoubi \cite[Theorem 2.3]{Kar2} for 
the so-called C-algebras. As was proved by H.~Leptin \cite{Lep} and 
J.~Wichmann \cite{Wic}, if $A$ is a symmetric Banach $\ast$-algebra, 
then the matrix algebra $M_{n}(A)$ is also symmetric. It follows that 
symmetric Banach $\ast$-algebras are C-algebras. The fact that 
$K_{n}(A) \simeq K_{n}(C^{\ast}(A))$ follows from the density theorem 
in $K$-theory. Indeed, since $A$ is symmetric, $A$ is dense and 
spectral invariant subalgebra of $C^{\ast}(A)$. Indeed, 
$a \in A$ has the same spectrum in $A$ or in $C^{\ast}(A)$.
\end{proof1}

For symmetric Banach $\ast$-algebras of finite stable rank we can 
give a better description of the Witt groups of $A$ in terms of 
$C^{\ast}(A)$. We recall that $\bar P_{2n}(A)$ is the quotient 
of $P_{2n}(A)$, the subset of indempotents
of $M_{2n}(A)$, given by the equivalence relation
$p_1 \approx p_2 \Leftrightarrow \exists \alpha\in GL_{2n}(A), 
\alpha p_1\alpha^{-1} = p_2$.

\begin{cor}
	Let $A$ be a symmetric Banach $\ast$-algebra. Suppose $s = tsr(A)$ is 
	finite. Then 
	$$ W_{n}(A) \simeq \bar P_{2s}(C^{\ast}(A)) \mbox{ for } n \mbox{ 
	even }$$
	and 
	$$ W_{n}(A) \simeq GL_{s}(C^{\ast}(A))/GL_{s}(C^{\ast}(A))_{0} \mbox{ for } n \mbox{ 
	odd. }$$
\end{cor}

\begin{proof1}
Recall $A$ is a dense and spectral invariant subalgebra of
$C^{\ast}(A)$. Then we have (cf. for instance \cite{Bad}) 
$$Bsr (A) \leq Bsr (C^{\ast}(A)) = tsr (C^{\ast}(A)) \leq tsr 
(A) = s \; .$$

According to \cite{Cor86} and \cite{Rie87}, we have $K_{0}(C^{\ast}(A)) 
\simeq P_{2s}(C^{\ast}(A))$
and 
$$K_{1}(C^{\ast}(A)) \simeq 
GL_{s}(C^{\ast}(A))/GL_{s}(C^{\ast}(A))_{0}\; .$$ 
The result now follows from the above Theorem and Bott periodicity 
theorem.
\end{proof1}

\end{document}